\documentclass{fundam}
\usepackage{url} 
\usepackage[ruled,lined]{algorithm2e}
\usepackage{graphicx}
\usepackage{amsmath,graphicx}

\allowdisplaybreaks[4]

\def\ml{l\kern-0.45ex\raise0.1ex\hbox{'}\kern-0.10ex}
\newtheorem{conject}[definition]{Conjecture}
\newtheorem{problem}[definition]{Problem}



\begin{document}

\setcounter{page}{99}
\publyear{22}
\papernumber{2105}
\volume{185}
\issue{2}

  \finalVersionForIOS

\title{On Local Antimagic Vertex Coloring for Complete Full $t$-ary Trees}


\author{Martin Ba\v{c}a, Andrea Semani\v{c}ov\'{a}-Fe\v{n}ov\v{c}\'{i}kov\'{a}
 \\
Department of Applied Mathematics and Informatics \\
Technical University,  Ko\v{s}ice, Slovak Republic, Slovakia\\
\{martin.baca, andrea.fenovcikova\}@tuke.sk
\and
Ruei-Ting Lai, Tao-Ming Wang\thanks{Address for correspondence:
           Department of Applied Mathematics, Tunghai University, Taichung, Taiwan, ROC. \newline \newline
          \vspace*{-6mm}{\scriptsize{Received June 2020; \ accepted April 2022.
 }}}
\\
Department of Applied Mathematics\\
Tunghai University, Taichung, Taiwan, ROC\\
s05240046@thu.edu.tw, wang@go.thu.edu.tw
}

 \maketitle

\runninghead{M. Ba\v{c}a  et al.}{Local Antimagic Vertex Coloring}

\begin{abstract}
Let $G = (V, E)$ be a finite simple undirected graph without $K_2$ components. A bijection $f : E \rightarrow
\{1, 2,\cdots, |E|\}$ is called a local antimagic labeling if for any two adjacent vertices $u$ and $v$, they have different vertex sums, i.e., $w(u) \neq w(v)$, where the vertex sum $w(u) = \sum_{e \in E(u)} f(e)$, and $E(u)$ is the set of edges incident to $u$. Thus any local antimagic labeling induces a proper vertex coloring of $G$ where the vertex $v$ is assigned the color (vertex sum) $w(v)$. The local antimagic chromatic number $\chi_{la}(G)$ is the minimum number of colors taken over all colorings induced by local antimagic labelings of $G$. It was conjectured \cite{Aru-Wang} that for every tree $T$ the local antimagic chromatic number $l+ 1 \leq \chi_{la} ( T )\leq l+2$, where $l$ is the number of leaves of $T$. In this article we verify the above conjecture for complete full $t$-ary trees, for $t \geq 2$. A complete full $t$-ary tree is a rooted tree in which all nodes have exactly $t$ children except leaves and every leaf is of the same depth. In particular we obtain that the exact value for the local antimagic chromatic number of all complete full $t$-ary trees is $ l+1$ for odd $t$.
\end{abstract}



\section{Background and introduction}
By a graph $G = (V, E)$ we mean a finite simple undirected graph without $K_2$ components.
For graph theoretic terminology we refer to Chartrand and Lesniak \cite{CL}. Hartsfield and Ringel \cite{Ad} introduced the concept of antimagic labeling of a graph, and conjectured that every connected graph except $K_2$ admits such an antimagic labeling, which remains unsettled till today.

\begin{definition}
Let $G = (V, E)$ be a finite simple undirected graph without isolated vertices. Let $f : E(G) \to \{1, 2, \cdots , |E(G)|\}$ be a bijection. For each
vertex $u \in V(G)$, the weight $w(u) = \sum_{e \in E(u)} f(e)$, where $E(u)$ is the set of edges
incident to $u$. If $w(u) \neq w(v)$ for any two distinct vertices $u$ and $v \in V(G)$, then
$f$ is called an {\bf antimagic} labeling of $G$. A graph $G$ is called antimagic if $G$ has an
antimagic labeling. Furthermore $f$ is called a {\bf local antimagic} labeling if $w(u) \neq w(v)$   for any two {\bf adjacent} vertices $u$ and $v \in V(G)$.
\end{definition}

Thus it is clear that any  antimagic labeling must be local antimagic. Moreover, a local antimagic labeling induces a proper vertex coloring of $G$ where the vertex $v$ is assigned the color (vertex sum) $w(v)$. This notation {\bf local antimagic labeling} was raised in 2017 by the following two sets of authors independently: Arumugam, Premalatha, Ba\v{c}a, Semani\v{c}ov\'a-Fe\v{n}ov\v{c}\'ikov\'a \cite{A}, and Bensmail, Senhaji and Lyngsie \cite{B}. Both groups also raised the following conjecture: Every connected graph other than $K_2$ is local antimagic. Bensmail {\it et al.} \cite{B} propose the slightly stronger form of the previous conjecture that every graph without component isomorphic to $K_2$ has a local antimagic labeling and proved that this conjecture is true for trees. However this conjecture has been proved by Haslegrave, using the probabilistic method, in \cite{H} more recently.

The {\bf local antimagic chromatic number} $\chi_{la}(G)$ of a graph $G$ is the minimum number of colors taken over all colorings induced by local antimagic labelings of $G$. Let $\chi(G)$ be the usual chromatic number of a graph $G$.
For any graph $G$, $\chi_{la}(G) \ge \chi(G)$. It was noted by Arumugam et al. \cite{A} that the difference $\chi_{la}(G) - \chi(G)$ can be arbitrarily large as shown in the following theorem.

\begin{theorem}\label{leaf-lower-bound}{\bf (\cite{A}, 2017)}
For any tree $T$ with $l$ leaves, $\chi_{la}(T) \ge l+1$.
\end{theorem}

Among others in \cite{A} the local antimagic chromatic number of paths, cycles, friendship graphs, wheels, complete bipartite graphs were studied. In this article we study the local antimagic vertex coloring of the trees. The local antimagic chromatic numbers for various corona products of $G$ with $\overline{K_m}$, when the graph $G$ is $P_n, C_n, K_n$, were reported in \cite{Aru-Wang}. It was furthermore conjectured in \cite{Aru-Wang} that for any tree $T$ with $l$ leaves, the local antimagic chromatic number $\chi_{la}(T)$ is either $l+1$ or $l+2$. In the following sections we verify the above conjecture for complete full $t$-ary trees, $t\ge 2$.  The definition of complete full $t$-ary trees is given below.

\begin{definition}
Let $t\geq 2$ be a positive integer. A {\bf complete full $t$-ary tree} is a rooted tree in which all nodes have exactly $t$ children except leaves and every leaf is of the same depth.
\end{definition}

\begin{figure}[ht!]
\centering
\includegraphics[scale=0.9]{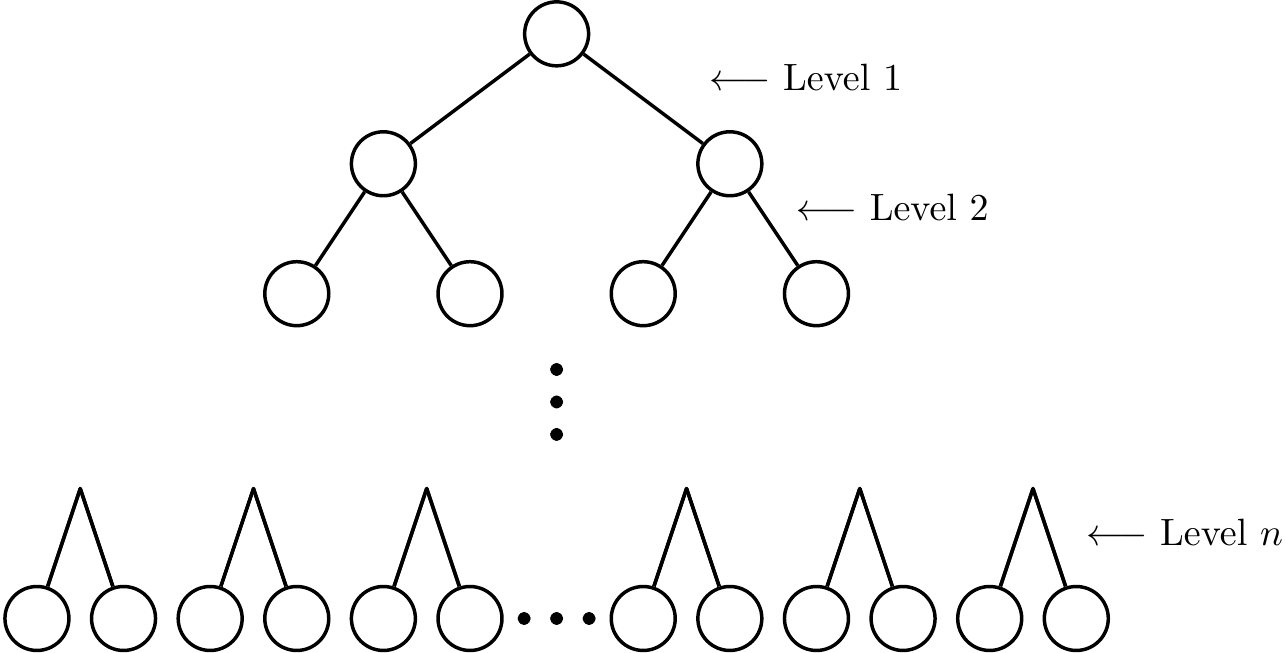}
\caption{Levels in a complete full binary tree.}\label{full-binary-tree}\vspace*{-3mm}
\end{figure}

\begin{definition}
The edge set of a complete full $t$-ary tree, $t \geq 2$, can be partitioned into $n$ levels, where the edges of level $l$, $1 \le l \le n,$ are the $l$-th edges from the root vertex. Note that there are $t^l$ such edges of level $l$.
\end{definition}
See Figure~\ref{full-binary-tree} for a binary case.

\section{Key observations}

In literature  Skolem \cite{Sk} remarked that Steiner triple systems
could be constructed from a sequence of integers $1,2,\cdots, 2n$ if these integers
could be arranged in disjoint pairs ($n$ of them) such that the differences are $1,2,\cdots, n$. He showed that such sequence exists if and only if $n \equiv 0,1\pmod 4$. For the purpose of arrangement for labeling, we study an additive analogue and extensions for this sequence as following. Note that we denote the set of consecutive integers $\{ m, m+1, \cdots, n \}$ by $[m,n]$.

\begin{problem}
For the first $2n$ positive integers $[1,2n]$, can one partition them into $n$ pairs so that the corresponding $n$ pair sums are consecutive?
\end{problem}
Note that the answer is positive if $n$ is odd and is negative if $n$ is even.  Let us see the case $n$ is odd first. Assume $n=2k+1$. Then we have the following.

\begin{lemma}\label{consective-sum-pair_4k+2}
Let $[1,4k+2]$ be the set of consecutive $4k+2$ integers, where $k$ is a positive integer. Then there exists a partition for $[1,4k+2]$ into $2k+1$ pairs for which the sums of these pairs form a set of $2k+1$ consecutive integers $[3k+3, 5k+3]$.
\end{lemma}
\begin{proof}
We consider the partition for $[1,4k+2]$ into $2k+1$ pairs as follows: $(1,3k+2), (3,3k+1), \cdots, (2k+1,2k+2)$ and $(2,4k+2), (4, 4k+1), (6,4k), \cdots, (2k, 3k+3)$. Therefore the set of pair sums are $[3k+3,  4k+3]\cup[ 4k+4,   5k+3]$, which is a set of $2k+1$ consecutive integers $[3k+3, 5k+3]$. This completes the proof. 
\end{proof}
The case when $n$ is even is related our study and we put the details in the following lemma.

\begin{lemma}\label{consective-sum-pair_4k}
Let $[1,4k]$ be the set of consecutive $4k$ integers, where $k$ is a positive integer. Then there exists a partition for $[1,4k]$ into $2k$ pairs for which the sums of these pairs form a set   of $2k-1$ consecutive integers $[3k+1,5k-1]$ with one isolated number $ 6k $. Moreover this result is optimal in the sense that for $[1,4k]$ there does not exist any partition into $2k$ pairs with consecutive pair sums.
\end{lemma}

\begin{proof}
We use the following partition for $[1,4k]$ into $2k$ pairs: $(1,3k), (3,3k-1), \cdots, (2k-1,2k+1)$ and $(2,4k-1), (4, 4k-2), (6,4k-3), \cdots, (2k-2, 3k+1)$. The remaining two integers are $(2k,4k)$. Therefore the set of pair sums are $[3k+1, 4k] \cup[4k+1,   5k-1] \cup \{ 6k \}$, which is a set of $2k-1$ consecutive integers with one isolated number $ 6k $.

As to see the result is optimal, one can conclude from the parity check as follows. Suppose $[1,4k]$ can be partitioned into pairs with consecutive sums, say $a, a+1, \cdots, a+(2k-1)$. Then $1+2+\cdots + 4k = a+ (a+1)+ \cdots+ (a+2k-1)$, and it follows that $2(4k+1) = 2a+2k-1$, which leads to a contradiction by parity. This completes the proof.
\end{proof}

For the above lemmas, it is not hard to see they may be extended by translation to the following general situations respectively for consecutive $4k$ and $4k+2$ integers to have consecutive pair sums.

\begin{lemma}\label{consective-sum-pair}
Let $[m,n]$ be the set of consecutive $n-m+1 = 4k$ integers, where $k$ be a positive integer. Then there exists a partition for $[m,n]$ into $2k$ pairs so that the pair sums form a set $[ m+n-k, m+n +k-2] \cup \{ m+n+2k-1 \}$, i.e., a set of $2k-1$ consecutive integers with one isolated number. Moreover this result is optimal in the sense that there does not exist any partition into $2k$ pairs with consecutive pair sums.
\end{lemma}

\begin{lemma}\label{consective-sum-pair-more}
Let $[m,n]$ be the set of consecutive $n-m+1 = 4k+2$ integers, where $k$ be a positive integer. Then there exists a partition for $[m,n]$ into $2k+1$ pairs for which the sums of these pairs form a set of $2k+1$ consecutive integers $[m+n-k, m+n+k]$.
\end{lemma}

More generally, one may ask the following.

\begin{problem}
For $t \ge 3$ and for first $tn$ positive integers $1,2,\cdots, tn$, can one partition them into $t$-tuples ($n$ of them) so that the sums of $t$-tuples are $n$ consecutive integers?
\end{problem}

For instance, for first $3n$ positive integers $1,2,\cdots, 3n$, can one partition them into $n$ triples so that the $n$ triple sums are consecutive? We have the following.

\begin{lemma}\label{consective-sum-triple}
For any $3n$ consecutive integers, $n\ge 1$, one can partition into $n$ triples so that the sums of the triples are $n$ consecutive integers.
\end{lemma}
\begin{proof}
Without loss of generality may assume the consecutive $3n$ in\-te\-gers to be $[1,3n]$. Write $[1,3n] = [1,n] \cup [n+1, 2n] \cup [2n+1, 3n] = A\cup B \cup C$, where $[1,n]  = A$, $[n+1, 2n]  = B $, and $ [2n+1, 3n] =  C$, respectively. The philosophy of the proof is based upon adding corresponding elements in $A$ and $B$ with reversed order to obtain constant partial sums, then put elements of $C$ to make the total sums of triples consecutive. Then the following partition into triples $(a,b,c)$, where $a\in A, b\in B, c\in C$, can do the job: $(1, 2n, 2n+1), (2, 2n-1, 2n+2), (3, 2n-2, 2n+3), \cdots, (n, n+1, 3n)$. Note that the set of the $n$ consecutive integers are $[4n+2, 5n+1]$. 
\end{proof}
Similarly it can be seen that the above lemma can be extended to the following.
\begin{lemma}\label{consective-sum-general-odd-tuple}
For any $tn$ consecutive integers, $t\geq 3$ odd, one can partition into $n$ of $t$-tuples so that the sums of the $t$-tuples are $n$ consecutive integers.
\end{lemma}
\begin{proof}Without loss of generality may assume the consecutive $tn$ in\-te\-gers to be $[1,tn]$.
 We partition $[1,tn]$ into $t$ consecutive subsets $S_1 \cup S_2 \cup \cdots \cup S_t$, where $S_1 = [1, n]$, $S_2 = [n+1, 2n]$, up to $S_t = [(t-1)n+1, tn]$. Then one may pick first $t-1$ elements of $t$-tuples in order from the pairs of subsets $(S_1, S_{2})$, $(S_3, S_{4})$, and so forth, up to $(S_{t-2}, S_{t-1})$. By the same philosophy while ob\-tai\-ning constant partial sums via adding reversed order elements and then pick the last element of the $t$-tuple from $S_t$ sequentially, one can see clearly the $n$ of the sums for each of these $t$-tuples are consecutive.
This completes the proof. 
\end{proof}

For even $t\geq 2$ and odd $n$, we have the following.
\begin{lemma}\label{consective-sum-general-even-tuple}
For any $tn$ consecutive integers, if $t \geq 2$ is even and $n$ is odd, one can partition into $t$-tuples so that these $n$ sums for each of the $t$-tuples are $n$ consecutive integers.
\end{lemma}
\begin{proof}Again without loss of generality we may assume the consecutive $tn$ in\-te\-gers to be $[1,tn]$.
If $t=2$, then by Lemma~\ref{consective-sum-pair_4k+2} we are done. If $t \geq 4$, then partition $[1,tn]$ into $t$ consecutive subsets $S_1 \cup S_2 \cup \cdots \cup S_t$, where $S_1 = [1, n]$, $S_2 = [n+1, 2n]$, up to $S_t = [(t-1)n+1, tn]$. Then one may pick $t$-tuples in order from subsets $(S_1, S_t)$, $(S_2, S_{t-1})$, and so forth to $(S_{{t}/{2}-1}, S_{{t}/{2}+2})$ in the following way. While obtaining constant partial sums $1+tn$ with adding reversed order elements and then apply Lemma~\ref{consective-sum-pair_4k+2} to the remaining last two $(S_{{t}/{2}}, S_{{t}/{2}+1})$, one can see clearly these $n$ sums for each of the $t$-tuples are $n$ consecutive integers.
This completes the proof. 
\end{proof}

For even $t\geq 2$ and even $n$, we have the following.
\begin{lemma}\label{consective-sum-general-even-4k-tuple}
For any $tn$ consecutive integers, if $t \geq 2$ is even and $n$ is also even, one can partition into $t$-tuples so that these $n$ sums for each of the $t$-tuples are $n-1$ consecutive integers with one isolated number. Moreover this result is optimal in the sense that there does not exist any partition into $t$-tuples with consecutive sums.
\end{lemma}
\begin{proof}Without loss of generality we may assume the consecutive $tn$ in\-te\-gers to be $[1,tn]$.
If $t=2$, then by Lemma~\ref{consective-sum-pair_4k} we are done. If $t \geq 4$, then partition $[1,tn]$ into $t$ consecutive subsets $S_1 \cup S_2 \cup \cdots \cup S_t$, where $S_1 = [1, n]$, $S_2 = [n+1, 2n]$, up to $S_t = [(t-1)n+1, tn]$. Then one can partition $[1,tn]$ into $t$-tuples in the following way. One may pick first $t-2$ coordinates of a $t$-tuple in reversed order from subsets $(S_1, S_t)$, $(S_2, S_{t-1})$, and so forth to $(S_{{t}/{2}-1}, S_{{t}/{2}+2})$, to obtain constant partial sums $1+tn$. Then apply Lemma~\ref{consective-sum-pair_4k} to the remaining last two $(S_{{t}/{2}}, S_{{t}/{2}+1})$ to get the last two coordinates of the $t$-tuple, one can partition $[1,tn]$ into $t$-tuples so that these $n$ sums for each of the $t$-tuples are $n-1$ consecutive integers with one isolated number.

As to see the result is optimal, one can conclude from the parity check as in Lemma~\ref{consective-sum-pair_4k}. This completes the proof. 
\end{proof}

To summarize with results of above lemmas, we have the following.
\begin{theorem}\label{consective-sum-general-thm}
For any $tn$ consecutive integers, one can partition them into $t$-tuples ($n$ of them) so that the sums of $t$-tuples are $n$ consecutive integers either for odd $t \geq 3$, or for even $t \geq 2$ and odd $n$. Moreover, if $t \geq 2$ is even and $n$ is also even, then for first $tn$ consecutive positive integers $1,2,\cdots, tn$ there does not exist any partition into $t$-tuples ($n$ of them) so that the sums of $t$-tuples are $n$ consecutive integers. However one can partition them into $t$-tuples ($n$ of them) so that the sums of $t$-tuples are $n-1$ consecutive integers with one isolated number.
\end{theorem}

\section{Labeling scheme}

For a complete full $t$-ary tree $T$, $t\geq 2$, we use the following labeling scheme to show the local antimagic chromatic number is $l+ 1 \leq \chi_{la} ( T )\leq l+2$, where $l$ is the number of leaves of $T$.

\medskip
Let $n$ be the number of levels of the complete full $t$-ary trees $T$. Note that if $n=1$, the local antimagic chromatic number is obviously $l+ 1$. In case $n=2$, it can be seen that the local antimagic chromatic number is at most $l+ 2$ since one may arrange the edge labels $1,2,\cdots, t$ for the level one and $t+1, t+2, \cdots, t+t^2$ in such a way that grouping $t^2$ into $t$ $t$-tuples with $t-1$ consecutive sums and one isolate by Lemma~\ref{consective-sum-pair_4k}, then adding them reversely with $1,2,\cdots, t$ to make $t-1$ constant sums and another isolate, which give rise to two new colors other than that on leaves. The top vertex receives the color $1+2+\cdots +t =  {t(t+1)}/{2}$ which is among $t+1, t+2, \cdots, t+t^2$ due to the inequality $t+1 \leq  {t(t+1)}/{2} \leq t+t^2$. Therefore the local antimagic chromatic number is either $l+1$ or $l+ 2$.

The labeling scheme for $n \ge 3$ and the sketch of the proof are given as follows.

\medskip
\noindent {\bf Step~1:} When $n$ is even, arrange labeling of
edges of levels by sequential order $2,4,\cdots, n-2, n-1, n-3, \cdots, 3,1,n$. Whenever $n$ is odd, arrange labeling of edges of levels by sequential order $2,4,\cdots, n-1, n-2, n-4, \cdots, 3,1,n$. Note that the edge labels are used sequentially by consecutive integers, and will be fixed in later steps by obtaining appropriate constant colors with possible isolated numbers.

\medskip
\noindent {\bf Step~2:} We fix the edge labels from the top level to the bottom level sequentially according to the requirement for colors. Start fixing the $t$ edge labels for the level 1 arbitrarily. Hence the color for the top vertex is done. Then in order to make constant colors (vertex sums) for vertices in the second row, we may use Lemma~\ref{consective-sum-general-odd-tuple} to arrange the $t^2$ edge labels in level 2. Inductively we can assign edge labels for levels one by one from the top to the bottom, so that we have the constant colors for vertices of the same row until the row right above the row of leaves.

 \medskip
\noindent {\bf Step~3:} It may be checked that the above labeling scheme is feasible, as the adjacent vertices receive the different color, the colors for all vertices till the third last row will be the colors of some leaves, and the vertices of the second last row will receive the constant color which is not among those colors for leaves. Therefore only one new color other than that of leaves and the result is obtained.

\medskip
More precisely the edge labeling is as follows. Let $n$ be the number of levels defined as before. When $n$ is even, arrange labeling of
edges of levels $2,4,\cdots, n-2, n-1, n-3, \cdots, 3,1,n$. Thus the range for integers used for these levels is $[1, t^2], [t^2 +1, t^2 + t^4], [t^2 + t^4 +1, t^2 + t^4 + t^6],$ up to $[1 + t^2 + t^3 +\cdots + t^{n-1}, t+ t^2 + t^3 +\cdots + t^n]= [({t^{n} - 1})/({t-1}), ({t^{n+1} - t})/({t-1})]$. We summarize up the edge labels level by level as following. A notion is used here for convenience.

\begin{definition}
Let $t\geq 2$ be a positive integer. For non-negative integers $n, k, a$, with $n \ge k \ge a \ge 0$ and $k, a$ are of the same parity, we denote by $(n,k,a)_t$ the abbreviation of the $t$-ary number $t^n + t^{n-1} + \cdots + t^{k+1} + t^k + t^{k-2} + t^{k-4} +\cdots + t^{a+4} + t^{a+2} + t^a$.
\end{definition}

\noindent{\bf Case~1:} When $n$ is even, arrange labeling of edges in the sequential order for levels $2,4,\cdots, n-2, n-1, n-3, \cdots, 3,1,n$.\smallskip\\
\noindent\begin{tabular}{ll}\small
Level 1:   &  $[(n-1,2,0)_t,(n-1,1,1)_t]$.\\
Level 2:   &  $[1,t^2]= [(0,0,0)_t, (2,2,2)_t]$.\\
Level 3:   &  $[(n-1,4,0)_t,(n-1,2,2)_t]$.\\
Level 4:   &  $[t^2 +1, t^2 + t^4] = [(2,2,0)_t, (4,4,2)_t]$.\\
Level 5:   &  $[(n-1,6,0)_t,(n-1,4,2)_t]$.\\
Level 6:   &  $[t^2 + t^4 +1, t^2 + t^4 + t^6] = [(4,4,0)_t,(6,6,2)_t]$.\\
Level 7:   &  $[(n-1,8,0)_t,(n-1,6,2)_t]$.\\
$\cdots$ & \\[-4pt]
$\cdots$ & \\
Level $n-4$: &    $[(n-6,n-6,0)_t,(n-4,n-4,2)_t]$.\\
Level $n-3$: &    $[(n-1,n-2,0)_t,(n-1,n-4,2)_t]$.\\
Level $n-2$:  &   $[(n-4,n-4,0)_t,(n-2,n-2,2)_t]$.\\
Level $n-1$:  &   $[(n-2,n-2,0)_t,(n-1,n-2,2)_t]$.\\
Level $n$:    &   $[\tfrac{t^{n} - 1}{t-1}, \tfrac{t^{n+1} - t}{t-1}] = [(n-1,0,0)_t,(n,1,1)_t]$.\smallskip\\
\end{tabular}

\medskip
\noindent{\bf Case~2:} Whenever $n$ is odd, arrange labeling of edges of levels $2,4,\cdots, n-1, n-2, n-4, \cdots, 3,1,n$. \\
\noindent\begin{tabular}{ll}\small
 & \\[-11pt]
Level 1:  &   $[(n-1,2,0)_t,(n-1,1,1)_t]$.\\
Level 2:  &   $[1,t^2]= [(0,0,0)_t, (2,2,2)_t]$.\\
Level 3:  &   $[(n-1,4,0)_t,(n-1,2,2)_t]$.\\
Level 4:  &   $[t^2 +1, t^2 + t^4] = [(2,2,0)_t, (4,4,2)_t]$.\\
Level 5:  &   $[(n-1,6,0)_t,(n-1,4,2)_t]$.\\
Level 6:  &   $[t^2 + t^4 +1, t^2 + t^4 + t^6] = [(4,4,0)_t,(6,6,2)_t]$.\\
Level 7:   &  $[(n-1,8,0)_t,(n-1,6,2)_t]$.\\
\end{tabular}

\noindent\begin{tabular}{ll}\small
$\cdots$& \\[-4pt]
$\cdots$& \\
Level $n-4$: &    $[(n-1,n-3,0)_t,(n-1,n-5,2)_t]$.\\
Level $n-3$: &    $[(n-5,n-5,0)_t,(n-3,n-3,2)_t]$.\\
Level $n-2$:  &   $[(n-1,n-1,0)_t,(n-1,n-3,2)_t]$.\\
Level $n-1$:  &   $[(n-3,n-3,0)_t,(n-1,n-1,2)_t]$.\\
Level $n$:    &   $[\tfrac{t^{n} - 1}{t-1}, \tfrac{t^{n+1} - t}{t-1}] = [(n-1,0,0)_t,(n,1,1)_t]$.\\
\end{tabular}

\begin{figure}[ht!]
\vspace{2mm}
\centering
\includegraphics[scale=0.9]{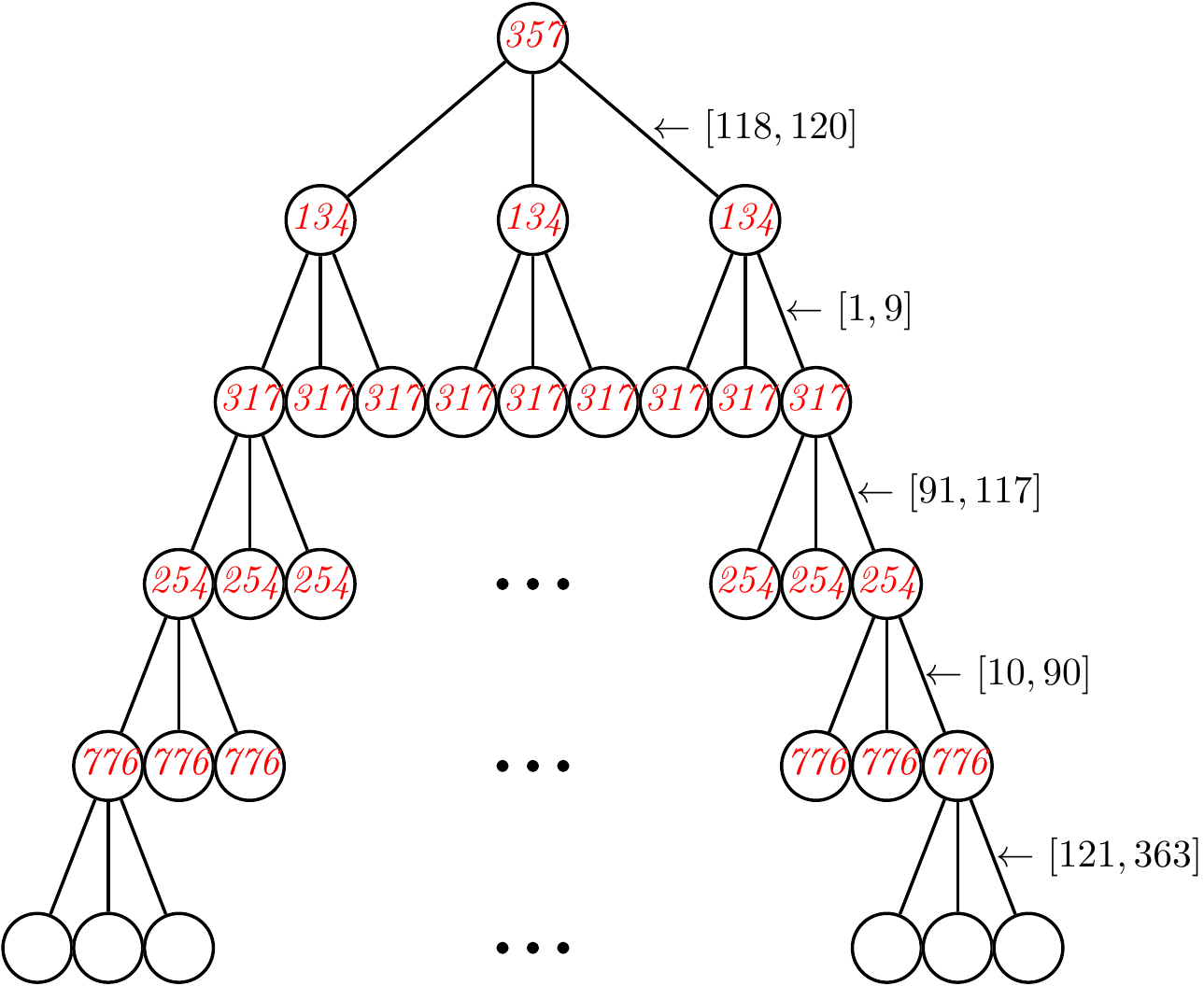}
\caption{Complete full $3$-ary tree with $5$ levels ($t=3$, $n=5$).}\label{full 3-ary tree of five levels}\vspace*{-4mm}
\end{figure}

\medskip
  Figure~\ref{full 3-ary tree of five levels} illustrates the labeling scheme and the range for edge labeling for the complete full 3-ary tree of five levels, that is the case $t=3$ and $n=5$. The precise way to get a proper coloring and associated chromatic number can be seen in the next section.

\section{Complete full $t$-ary trees for odd $t$}

In this section we calculate the local antimagic chromatic numbers for complete full $t$-ary trees, $t\geq 3$ odd, is shown to be $ l+1$, where $l$ is the number of leaves of $T$. In particular, for the class of complete full $t$-ary trees, where $t\geq 3$ odd, the conjecture that for a tree $T$ the local antimagic chromatic number $l+ 1 \leq \chi_{la} ( T )\leq l+2$ is verified as a byproduct.
\begin{theorem}\label{3-ary-main}
For complete full $t$-ary trees $T$, $t\geq 3$ odd, the local antimagic chromatic number $\chi_{la}(T) = l+1$, where $l$ is the number of leaves of $T$.
\end{theorem}

\begin{proof}
Let $n$ be the number of levels in $T$.
Follow the labeling scheme in previous section first. Then start arranging the labeling by fixing the $t$ edge labels for the level 1 arbitrarily. Hence the color (vertex sum) for the first row (the top vertex) is done. Then we may use Lemma~\ref{consective-sum-general-odd-tuple} to arrange the $t^2$ edge labels in level 2 as $t$ consecutive $t$-tuple-sums and add them to corresponding integers reversely in level 1, to make constant colors for vertices in the second row. Inductively we can assign edge labels for levels 3, 4, etc. one by one while applying Lemma~\ref{consective-sum-general-odd-tuple} again so that we obtain the constant vertex sums (colors) for vertices of the same row, until the row right above the row of leaves.

\medskip
Without loss of generality we use the case $n$ is even for the sketch of the proof. Let $h$ be odd and $1 \leq h \leq n-1$, then the calculation shows that the constant color for the vertices between the level $h-1$ and level $h$ is $( 2 t^{n+1} + (- t^3 + t^2 - t +1)\cdot t^{h-1} - t^2 - 1)/{(2t-2)}$ for $1 \leq h \leq n-1$. Note that in particular the vertex between level 0 and level 1 is the unique top vertex. Again the calculation shows that the constant color for the vertices between the level $h$ and level $h+1$ is $(2 t^{n} + (t^3 - t^2 + t -1)\cdot t^{h} - t^2 - 1)/{(2t-2)}$ for $1 \leq h \leq n-3$. We calculate the local antimagic chromatic number by the following steps.

\medskip
\noindent{\bf Step~1:} Show that the colors for the vertices of two adjacent rows are different.

\smallskip
Suppose on the contrary, for $1 \leq h \leq n-3$ we have $$ 2 t^{n+1} + (- t^3 + t^2 - t +1)\cdot t^{h-1} - t^2 - 1 =  2 t^{n} + (t^3 - t^2 + t -1)\cdot t^{h} - t^2 - 1$$ which implies $-(t^2 +1 )(\frac{t-1}{2})\cdot t^{h-1} = t^{n}$, a contradiction due to the parity. For $h = n-2$ (that is the row between level $n-1$ and $n-2$, and the row between level $n-2$ and $n-3$), one sees that $$ 2 t^{n} + (- t^3 + t^2 - t +1)\cdot t^{n-2} =  2 t^{n} + (t^3 - t^2 + t -1)\cdot t^{n-4}$$ which implies $(t-1 )\cdot (t^2(t^2 - 2) +1)\cdot  t^{n-4} = 0$, a contradiction since $t\geq 3$ odd.

\medskip
\noindent{\bf Step~2:} Show that the colors for non-leaves of all rows except the last two rows are among colors of leaves.

\smallskip
Notice that the colors for leaves are among the interval $[ ({t^n - 1})/({t-1}),$ $ ({t^{n+1} - t})/({t-1})]$. On the other hand, these constant colors for non-leaves (except the vertices over the second last rows) are among those of some leaves, since it is not hard to check the following inequalities hold: $$ \tfrac{t^n - 1}{t-1} \leq \tfrac{1}{2(t-1)}( 2 t^{n+1} + (- t^3 + t^2 - t +1)\cdot t^{h-1} - t^2 - 1) \leq \tfrac{t^{n+1} - t}{t-1}$$ for $1 \leq h \leq n-1$, and $$ \tfrac{t^n - 1}{t-1} \leq \tfrac{1}{2(t-1)}( 2 t^{n} + (t^3 - t^2 + t -1)\cdot t^{h} - t^2 - 1) \leq \tfrac{t^{n+1} - t}{t-1}$$ for $1 \leq h \leq n-3$, respectively.

\eject
\noindent{\bf Step~3:} Show that overall there is only one new color other than those colors of leaves.

\smallskip
We can see that the constant color for the vertices of the second last row is the only new color other than those on leaves. By Lemma~\ref{consective-sum-general-odd-tuple} it is obvious that this constant color is new, since it is the sum of some selected $t$-tuple from the bottom level $n$ with one corresponding edge of the level $n-1$, which is clearly strictly greater than $({t^{n+1} - t})/({t-1})$, since the largest edge label $({t^{n+1} - t})/({t-1})$ is among one element of the $t$-tuples selected to make the constant color. Therefore we are done with the proof for the result.
\end{proof}

See again Figure~\ref{full 3-ary tree of five levels} for  the example of the labeling range and calculated coloring (numbers in red) for the complete full 3-ary tree of five levels.

\section{Complete full $t$-ary trees for even $t$}

In this section, using the labeling scheme mentioned above, we verify the conjecture for the class of complete full $t$-ary trees for even $t\geq 2$.
\begin{theorem}\label{binary-tree-local-antimagic}
For complete full $t$-ary trees $T$, $t\geq 2$ even, the local antimagic chromatic number $\chi_{la}(T)$ is either $l+1$ or $l+2$, where $l$ is the number of leaves of $T$.
\end{theorem}

\begin{proof}
Let $n$ be the number of levels in $T$.
It suffices to show that $\chi_{la}(T) \le l+2$ for complete full $t$-ary trees $T$, $t\geq 2$ even, since we have the lower bound $l+1$ from Theorem~\ref{leaf-lower-bound} in \cite{A}.
Without loss of generality we use the case $n$ is even for the sketch of the proof.

\medskip
Let $h$ be odd and $1 \leq h \leq n-1$, then the calculation shows the constant {\bf non-isolated color N} (resulted from adding the consecutive $t$-tuple-sums and reversely consecutive numbers) and the isolated {\bf jump color J} (resulted from adding the smallest consecutive $t$-tuple-sums and reversely the largest of consecutive numbers) for the vertices between the level $h-1$ and level $h$ are $ ( 2 t^{n+1} + (- t^3 + t^2 - t +1)\cdot t^{h-1} - t^2 - 3t +2)/{(2t-2)}$ and $ ( 2 t^{n+1} + ( -t^3 + t^2 +2t -2)\cdot t^{h-1} - t^2 - 3t + 2)/{(2t-2)}$ respectively for $1 \leq h \leq n-1$.

Again the calculation shows that the constant non-isolated color and the isolated color for the vertices between the level $h$ and level $h+1$ are $( 2 t^{n} + (t^3 - t^2 + t -1)\cdot t^{h} - t^2 - 3t + 2)/{(2t-2)}$ and $ ( 2 t^{n} + (t^3 - t^2 + 4t -4)\cdot t^{h} - t^2 - 3t + 2)/{(2t-2)}$ respectively for $1 \leq h \leq n-3$.

As for the boundary row of vertices between level $n-1$ and $n-2$, the calculation shows that the constant non-isolated color and the isolated color are $ ( 2 t^{n} + (t^3 - t^2 + t -1)\cdot t^{n-2} - t^2 - t)/{(2t-2)}$ and $ (  t^{n+1} + (t^2 + 2t -2)\cdot t^{n-2} - t^2 - 3t + 2)/{(2t-2)}$ respectively.

\medskip
See Figure~\ref{full-binary-n=4} and Figure~\ref{full-binary-n=5} for the proof of the theorem as examples for binary cases, $n=4$ and $n=5$, respectively.

\begin{figure}[!h]
\vspace*{1mm}
\centering
\includegraphics[scale=0.9]{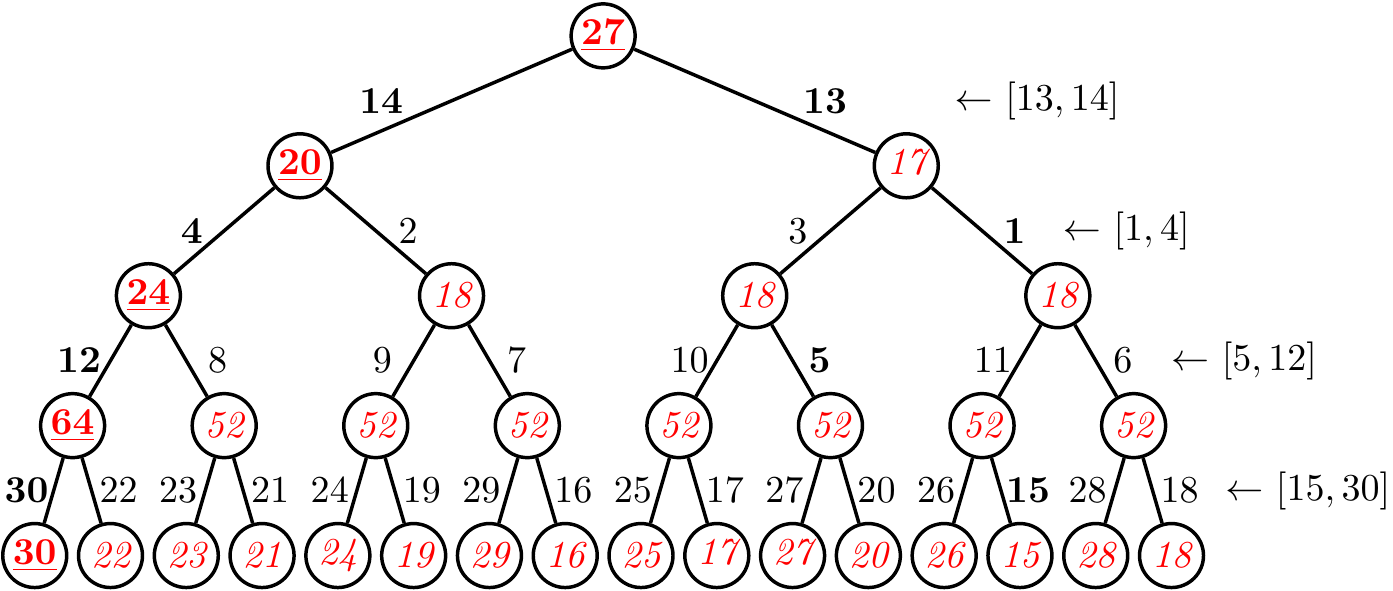}\vspace{1mm}
\caption{Complete full binary tree of four levels ($t=2$, $n=4$).}\label{full-binary-n=4}
\end{figure}
\begin{figure}[!h]
\vspace*{1mm}
\centering
\includegraphics[scale=0.9]{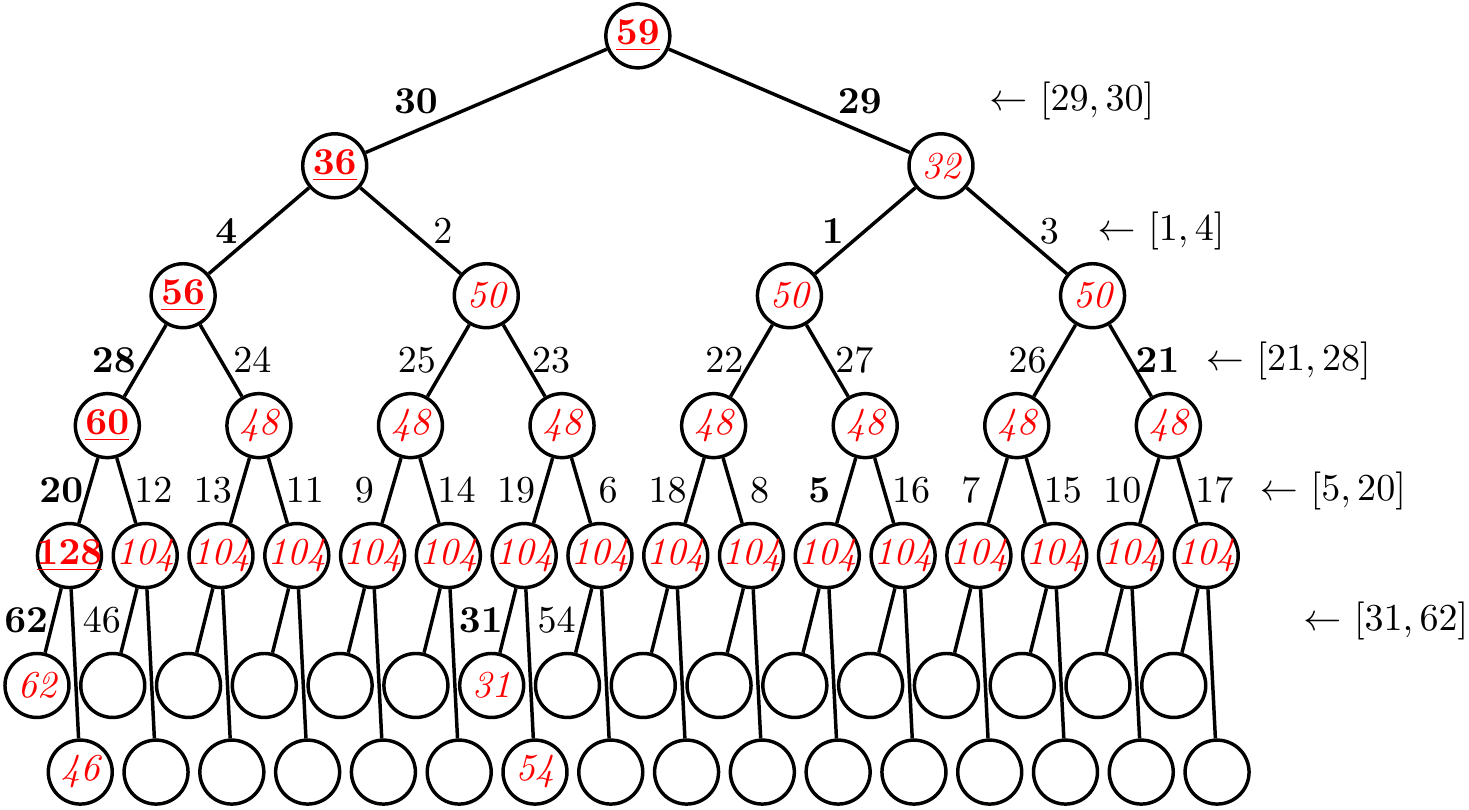}
\caption{Complete full binary tree of five levels ($t=2$, $n=5$).}\label{full-binary-n=5}
\end{figure}

\smallskip
We calculate the local antimagic chromatic number by the following steps.

\medskip
\noindent{\bf Step~1:} Show that the colors for the vertices of two adjacent rows are different.

\smallskip
Then we may check that the colors for the vertices of two adjacent rows are different. We split the situation into six cases as follows according to all possibilities among the constant non-isolated colors $N$ and the isolated colors $J$ for two adjacent rows (from vertices of the $h$-th row to that of the next lower $(h+1)$-th row), $1 \leq h \leq n-2$. Note that $n$ is even.

\medskip
\noindent{\bf Case~1: from $J$ to $J$ for $h \leq n-2$.}

\smallskip
 Suppose   the isolated color $J$ and the isolated color $J$ are the same. In this case, simply switch the $t$-tuple in the lower $(h+1)$-th level connecting the isolated color $J$ in the $h$-th level with the $t$-tuple in the lower $(h+1)$-th level connecting the smallest number in the $h$-th level. All the other edge labels stay the same. Then this will resolve the conflict of the case.

\eject 
\noindent{\bf Case~2: from $J$ to $N$ for $h \le n-3$.}

\smallskip
 Suppose the isolated color $J$ and the non-isolated color $N$ are the same. In this case, one has $$  2 t^{n+1} + ( -t^3 + t^2 +2t -2)\cdot t^{h-1} - t^2 - 3t + 2 =   2 t^{n} + (t^3 - t^2 + t -1)\cdot t^{h} - t^2 - 3t + 2 $$ which can be simplified to $$ 2t^{n-h+1} = t^3 + t^2 + t - 2.$$ This leads a contradiction due to the following. Since ${n-h+1} \geq 4$, we have $ 2t^{n-h+1} \geq 2t^4 $. Let $f(t)= 2t^4 - (t^3 + t^2 + t - 2)$. Note that for $t \geq 2$ one has $f'(t)=8t^3 - 3t^2 - 2t - 1 > 0$, and $f(2)=20$, thus $f(t)>0$, namely $2t^4 > t^3 + t^2 + t - 2 $. Therefore $ 2t^{n-h+1} \geq 2t^4 > t^3 + t^2 + t - 2 $, a contradiction.

\medskip
\noindent{\bf Case~3: from $N$ to $J$ for $h \leq n-3$.}

\smallskip
 Suppose   the non-isolated color $N$ and the isolated color $J$ are the same. In this case, one has $$ 2 t^{n+1} + (- t^3 + t^2 - t +1)\cdot t^{h-1} - t^2 - 3t +2  =   2 t^{n} + (t^3 - t^2 + 4t -4)\cdot t^{h} - t^2 - 3t + 2 $$ which can be simplified to $$ 2t^{n-h+1}  = t^3 + t^2 + 4t -1,$$ again this reaches a contradiction due to the parity.

\medskip
\noindent{\bf Case~4: from $N$ to $N$ for $h \leq n-2$.}

\smallskip
 Suppose   the non-isolated color $N$ and the non-isolated color $N$ are the same. For $1 \leq h \leq n-3$ we have $$  2 t^{n+1} + (- t^3 + t^2 - t +1)\cdot t^{h-1} - t^2 - 3t +2 =  2 t^{n} + (t^3 - t^2 + t -1)\cdot t^{h} - t^2 - 3t + 2$$ which implies $ 2t^{n-h+1} = t^3 + t^2 + t +1,$ a contradiction due to the parity.

\medskip
\noindent{\bf Case~5: from $N$ to $J$ for $h = n-2$.}

\smallskip
 This is for the row between level $n-2$ and $n-3$, and the row between level $n-1$ and $n-2$. Suppose   the non-isolated color $N$ and the isolated color $J$ are the same. One sees that $$ t^{n+1} + (t^2 + 2t -2)\cdot t^{n-2} - t^2 - 3t +2 =  2 t^{n} + (t^3 - t^2 + t -1)\cdot t^{n-4}- t^2 - t  $$ which can be simplified to $$(t^4 + t^2 -1 )\cdot t^{n-4} = 2,$$ a contradiction since $t^4 + t^2 -1$ is at least 19 and $n-4 \geq 0$.

\medskip
\noindent{\bf Case~6: from $J$ to $N$ for $h = n-2$.}

\smallskip
 Again this is for the row between level $n-2$ and $n-3$, and the row between level $n-1$ and $n-2$. Suppose   the isolated color $J$ and the non-isolated color $N$ are the same. One sees that $$ 2 t^{n} + (t^3 - t^2 - t +1)\cdot t^{n-2} - t^2 - t  =  2 t^{n} + (t^3 - t^2 + 4t -1)\cdot t^{n-3}- t^2 - 3t +2 $$ which can be simplified to $$ -(t^4 - 2t^3 - 3t +1)\tfrac{t^{n-3}}{2} = t-1,$$ this reaches a contradiction due to the parity if $n > 4, t \geq 2$. Also it is impossible for the remaining situation $n=4, t=2$.

\medskip
\noindent{\bf Step~2:} Show that the colors for non-leaves of all rows except over last two rows are among colors of leaves.

\smallskip
Notice that the colors for leaves are among the interval $[ ({t^n - 1})/({t-1}),$ $({t^{n+1} - t})/({t-1}) ]$. On the other hand, these constant colors for non-leaves (except the vertices over the second last rows) are among those of some leaves, since it is not hard to check that the following inequalities hold: $$ \tfrac{t^n - 1}{t-1} \leq \tfrac{1}{2(t-1)}( 2 t^{n+1} + (t^2 - t^3 - t +1)\cdot t^{h-1} - t^2 - 1) \leq \tfrac{t^{n+1} - t}{t-1}$$ for $1 \leq h \leq n-1$, and $$ \tfrac{t^n - 1}{t-1} \leq \tfrac{1}{2(t-1)}( 2 t^{n} + (t^3 - t^2 + t -1)\cdot t^{h-1} - t^2 - 1) \leq \tfrac{t^{n+1} - t}{t-1}$$ for $1 \leq h \leq n-3$, respectively.

\medskip
\noindent{\bf Step~3:} Show that there are only two new colors other than those colors of leaves.

\smallskip
Lastly we can see that, for the vertices of the second last row, there are only two new colors other than those colors of leaves. By Lemma~\ref{consective-sum-general-even-4k-tuple} one may arrange edge labels of the last level as $t^{n-1}$ pairs so that they form $t^{n-1} - 1$ consecutive integers and one isolated number. Then adding these integers to the consecutive edge labels for the $(n-1)$-th level in a reverse way as usual to get two constants, will give two colors of the   vertices of the second last row. Note that both colors are again strictly greater than the largest edge label $({t^{n+1} - t})/({t-1})$ as before. Therefore only two new colors other than those colors of leaves will be created and we are done with the proof.
\end{proof}

\begin{figure}[ht!]
\centering
\includegraphics{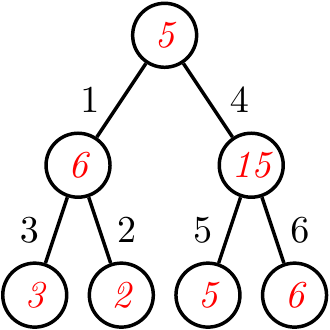}
\caption{Complete full binary tree of two levels ($t=2$, $n=2$).}\label{full-binary-n=2}
\end{figure}

\begin{figure}[ht!]
\centering
\includegraphics{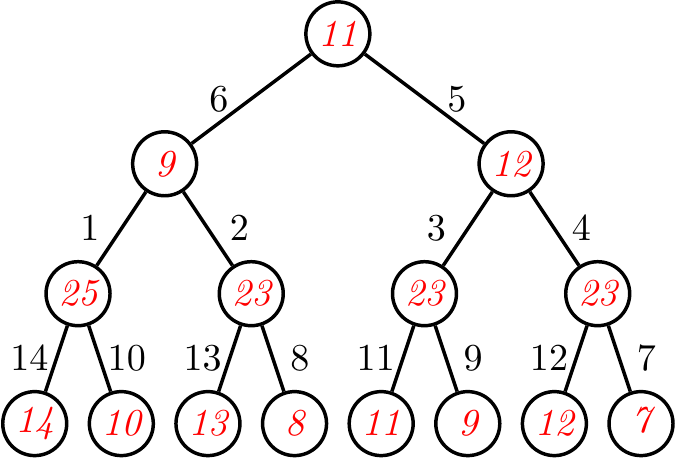}
\caption{Complete full binary tree of three levels ($t=2$, $n=3$).}\label{full-binary-n=3}
\end{figure}

The feasible labeling and coloring for $n=2,3$ are in the Figure~\ref{full-binary-n=2} and Figure~\ref{full-binary-n=3}, respectively. Note that for $n=2$ the coloring we obtain is of type $l+1$. In fact, while $l+2$ colors are used in the proof, however we believe that the exact value of the local antimagic chromatic number for the complete full $t$-ary tree, $t \geq 2$ even, is $l+1$, where $l$ is the number of leaves. See {Conjecture~\ref{even-t-conjecture}} in the next concluding section.


\section{Future studies}

In this article we are concerned with the local antimagic coloring of trees. It is interesting to observe that, one may use the properties regarding the additive version of the Skolem sequence to give rise to the arrangement of the labeling. Again this shows the intimate relationship between the integer sequences and graph labeling and coloring.

Among others, we obtain the upper bound $l+2$ for the local antimagic chromatic numbers of a complete full $t$-ary tree $T$ for all $t\geq 2$, and in particular we obtain that the exact value for the local antimagic chromatic number of a complete full $t$-ary tree $T$ is $ l+1$ for odd $t$, where $l$ is the number of leaves of $T$.
However the following conjectures are still open and need further exploration.

\begin{conject}\label{even-t-conjecture}
For any complete full $t$-ary tree $T$, $t\geq 2$ even,  $\chi_{la} ( T ) = l+1$, where $l$ is the number of degree one vertices (leaves) of $T$.
\end{conject}

\begin{conject}\label{previous-conjecture}
For any tree $T$ on at least three vertices,  $l+ 1 \leq \chi_{la} ( T )\leq l+2$, where $l$ is the number of degree one vertices (leaves) of $T$.
\end{conject}

Note that the Conjecture~\ref{previous-conjecture} was mentioned in \cite{Aru-Wang}.
It will be interesting to characterize trees $T$ such that $\chi_{la} ( T )= l+1$, where $l$ is the number of degree one vertices (leaves) of $T$.

\subsection*{Acknowledgement}

The research of  Martin Ba\v{c}a  and  Andrea Semani\v{c}ov\'{a}-Fe\v{n}ov\v{c}\'{i}kov\'{a}
 for this article was supported by APVV-19-0153 and by VEGA 1/0233/18.

\smallskip
The research for the author Tao-Ming Wang is supported by MOST 108-2918-I-029-002 and MOST 109-2115-M-029-005-MY2 from the ministry of science and technology of the government of Taiwan,  and this work was done  during
his visit to the Department of Applied Mathematics and Informatics, Technical University, Ko\v{s}ice, Slovakia.
Tao-Ming Wang wishes to express his sincere gratitude for the hospitality from the department and local hosts.

\smallskip
The authors appreciate very much the efforts that referees have made for corrections of this paper.


\end{document}